\documentclass{amsart}

\usepackage{latexsym,amssymb,amsmath, amsfonts,epsf}
\usepackage{tabularx}
\usepackage{amsthm}
\usepackage{amscd}

\newcommand{\R}{\mathbb{R}}
\newcommand{\C}{\mathbb{C}}

\newcommand{\la}{\Lambda}

\newcommand{\lb}{\nu}

\newcommand{\ls}{\lesssim}
\newcommand{\gs}{\gtrsim}

\newcommand{\de}{\delta}
\newcommand{\M}{\mathcal{M}}
\newcommand{\PP}{\mathcal{P}}
\newcommand{\un}{\underline{\theta}}

\newtheorem{theorem}{Theorem}
\newtheorem{lemma}[theorem]{Lemma}

\address{Department of Mathematics, Universidad Aut\'onoma de Madrid, 28049 
Madrid, Spain}\email{K.M.Rogers.99@cantab.net}
\author{Keith M. Rogers}
\date{}
\title{On a planar variant of the Kakeya problem}
\thanks{Supported by the EC project HPRN-CT-2001-00273 - HARP}

\begin{document}

\begin{abstract}  
 A {\it $\mathcal{K}^n_2$-set} is a set of zero Lebesgue measure containing a 
translate of every plane in an $(n-2)$--dimensional manifold in $\mathrm{Gr}(n,2)$, 
where the manifold fulfills a curvature condition.      
We show that this is a natural 
class of sets with respect to the Kakeya problem and prove that $\dim_H(E)\ge 
7/2$ for all $\mathcal{K}^4_2$-sets~$E$.  When the 
underlying field is replaced by $\C$, we get $\dim_H(E)\ge 7$ 
for all $\mathcal{K}^4_2$-sets over $\C$, 
and we construct an example to show that this is sharp. Thus 
$\mathcal{K}^4_2$-sets over $\C$ do not necessarily have full Hausdorff 
dimension.   \end{abstract}

\maketitle

\section{Introduction}\label{intro}

A {\it Besicovitch set (with lines)} in $\R^n$ is a set of zero Lebesgue measure 
that contains a translate of every member of $\mathrm{Gr}(n,1)$, 
where $\mathrm{Gr}(n,1)$ is the Grassmanian manifold of 1--dimensional linear 
subspaces of~$\R^n.$ 

A construction of  A. Besicovitch \cite{be} led to the surprising fact that such sets 
exist. The Kakeya
conjecture asserts that Besicovitch sets must have full Hausdorff dimension.

We note that $\mathrm{Gr}(n,1)$ has dimension~$n-1$ and a line has dimension~1, 
and the union of these lines will fill the space and have dimension~$n$. 
Informally, the Kakeya conjecture asserts that under translation of the lines, the dimensions 
continue to add up, or that the intersection remains 
negligible.

The conjecture has been solved in the affirmative in the plane, but 
is open for $n\ge3$.
T.~Wolff~\cite{wo} proved that $\dim_H(E) \ge(n+2)/2$ for all Besicovitch sets 
$E$ in $\R^n,$ where  $\dim_H$ denotes Hausdorff dimension, 
and there has been much recent progress for higher 
dimensions (see \cite{bo2}, \cite{ka1}, \cite{ka2}, \cite{la}, \cite{ta}).

N. H. Katz, 
I. \L aba\ and\ T. Tao \cite{ka} have shown that the Minkowski dimension of a 
Kakeya set in $\R^3$ is strictly greater than $5/2$. They also show that if a 
Kakeya set in $\R^3$ has dimension close to $5/2$, then it must exhibit a certain 
structural property that they call `planiness'. Roughly speaking,
 most of the lines that pass through 
a point in the Besicovitch set, lie in the union of a small number of planes. 
It seems reasonable then, to consider the variant of the problem where lines 
are replaced by planes. 

An {\it $(n,2)$-set} is a set that contains a translate of every member of 
$\mathrm{Gr}(n,2)$, where $\mathrm{Gr}(n,2)$ is the Grassmanian manifold of 
2--dimensional linear subspaces of~$\R^n.$ 
J.~Bourgain~\cite{bo} proved that $(4,2)$-sets have strictly positive
Lebesgue measure, so that there are no, so called, Besicovitch $(4,2)$-sets.
Similarly, there are no Besicovitch $(3,2)$-sets.  
T.~Mitsis~\cite{mi} recently claimed that $(n,2)$-sets have full Hausdorff 
dimension when $n\ge5$, although unfortunately the argument is incomplete.

As $\mathrm{Gr}(n,2)$ has dimension $2(n-2)$ and a plane has dimension 2, these 
results do not have the same informal interpretation as that of the 
Kakeya problem. Indeed the planes of an $(n,2)$-set will inevitably have, 
in some sense, nontrivial intersection. 

Initially, it appears reasonable to ask  
whether sets containing translates of an \newline $(n-2)$--dimensional
manifold in $\mathrm{Gr}(n,2)$ necessarily have full Hausdorff dimension. 
Without further restriction on the manifold however,  
this will fail. For instance, a 2--dimensional manifold in $\mathrm{Gr}(4,2)$ 
could have all of its planes contained in a 3--space, so that the set consisting of 
the union of these planes would have Hausdorff dimension 3. We outlaw this by 
adding a curvature condition which is essentially a version of the Wolff axiom 
(\cite{wo}, \cite{ka}, \cite{mo}) for planes.

 Let $\pi,\pi'\in\mathrm{Gr}(n,2),$ and define the {\it major angle between $\pi$ and $\pi'$} by 
$$
\theta(\pi,\pi')=\|\mathrm{proj}_\pi-\mathrm{proj}_{\pi'}\|,
$$
where $\mathrm{proj}_\pi:\R^n\to\pi$ is the orthogonal projection onto $\pi$. We note that 
$\theta$ takes values from zero to one, so this is not a standard angle.  
As $\theta$ is a metric, we can induce a measure $\mu$ on the manifold $\mathcal{M}\subset\mathrm{Gr}(n,2)$, 
which we normalise to have total mass one. We denote by $B_\epsilon(\pi)$ 
the ball 
$$\{\pi'\in\mathrm{Gr}(n,2): \theta(\pi,\pi')<\epsilon\},$$ and denote  
by $\mathcal{M}^{\epsilon}_{\la_m}$ the set
$$\{\pi'\in\mathcal{M}:\min_{\pi\subset\la_m}\theta(\pi,\pi')<\epsilon\},$$ 
where $\la_m\subset\R^n$ 
is an $m$--space with $3\le m\le n-1$.

An $(n-2)$--dimensional  manifold 
$\mathcal{M} \subset \mathrm{Gr}(n,2)$ is said to be {\it curved} if 
there exists a constant $C$ such that
\begin{equation}\label{pope}
\frac{\mu(\mathcal{M}^{\epsilon_1}_{\la_m}\cap 
 B_{\epsilon_2}(\pi))}{\epsilon_1^{m-2}\epsilon_2^{n-m}}<C
\end{equation}
for all $\epsilon_1,\epsilon_2>0,$  
$m$--spaces $\la_m\subset\R^n$ and $\pi\in\mathcal{M}.$

The condition forces 
the  dimension of $\{\pi\in\mathcal{M}:\pi\subset\Lambda_m\}$ to be less than or 
equal to $m-2$ for
all $m$--spaces $\la_m,$ and
the manifold to be more evenly distributed in~$\mathrm{Gr}(n,2)$. 
The reason that this is known as an axiom in the Kakeya problem, 
is that the corresponding condition is  
automatically fulfilled by $\mathrm{Gr}(n,1)$. This is not 
the case for planes, so we are obliged to include it in our definition. 

A {\it $\mathcal{K}^n_2$-set} is a set of zero Lebesgue measure that contains a 
translate of every member of an $(n-2)$--dimensional curved manifold 
$\mathcal{M}\subset\mathrm{Gr}(n,2).$
 
 As $\{ l\times \R : l\in \mathrm{Gr}(n-1,1)\}$ is an $(n-2)$--dimensional curved manifold, 
 a Besicovitch set that is 
constructed by taking the cross product of a Besicovitch set in $\R^{n-1}$ with 
a copy of $\R$ is also a $\mathcal{K}^n_2$-set. Thus we have the existence of 
such sets.

In Section~\ref{two} we will prove the following results.
\begin{theorem}\label{dimension0}
Let $E$ be a $\mathcal{K}^3_2$-set. Then $\dim_H(E)=3.$
\end{theorem}

Thus $\mathcal{K}^3_2$-sets  have full Hausdorff dimension. Our main concern will be 
the proof of the following result. 

\begin{theorem}\label{dimension}
Let $E$ be a $\mathcal{K}^4_2$-set. Then $\dim_H(E)\ge 
7/2$.
\end{theorem}

If we replace the underlying field by $\C$ and multiply the exponents of $\epsilon_1$ and $\epsilon_2$ in (\ref{pope}) by two, 
then the proof of Theorem~\ref{dimension}
can be modified to obtain the following bound.   
The modification involves little more than changing the relevant exponents.

\begin{theorem}\label{cork}
Let $E$ be a $\mathcal{K}^4_2$-set over $\C$. Then 
$\dim_H(E)\ge 7$.
\end{theorem}

In Section~\ref{three}, we will construct an example to show that Theorem~\ref{cork} is sharp. Thus $\mathcal{K}^4_2$-sets over $\C$ do not necessarily have 
full Hausdorff dimension.

\section{Geometric preliminaries}\label{one}

Throughout $\de$ will be a real parameter such that $0<\de\ll 1$, and $\pi$ will 
denote a plane, and never the usual number. We say that a set  
$\Pi\subset\mathrm{Gr}(n,2)$ is {\it $\delta$--separated} if  
$\theta(\pi,\pi')> \delta$ for all $\pi,\pi'\in\Pi$.

 We use $A\ls B$ to denote the estimate $A\le C_\epsilon\de^{-\epsilon} B$ for 
all $\epsilon>0$, where $C_\epsilon$ is a constant depending only on $\epsilon$ and the manifold $\mathcal{M}$. 
This notation  will be convenient, as factors of $\log(1/\de)$ will simply 
disappear.  
We use $A\sim B$ to denote   $B/2< A\le B.$   

We will require another notion of an angle between two planes. 
The need for two notions is created by the fact that 
in $\R^n$, where $n\ge4$, the intersection of two planes can be a point, 
and not necessarily a line.
Let $l, l'\in \mathrm{Gr}(n,1),$ and define the {\it minor angle between $\pi$ and $\pi'$} 
 by
$$
 \underline{\theta}(\pi,\pi')=\min_{l\subset \pi,l'\subset 
\pi'}\|\mathrm{proj}_l-\mathrm{proj}_{l'}\|. $$
  Informally, $\underline{\theta}$ can be considered to be the smallest angle 
between two planes. If $\pi,\pi'\in\mathrm{Gr}(n,2)$ intersect in a line, then $\underline{\theta}(\pi,\pi')=0$, so that $\underline{\theta}$ is not a metric.
  
Define a {\it plate} $P_\pi$ to be the image of 
$[0,1]\times[0,1]\times\underbrace{[0,\de]\times\cdots\times[0,\de]}_{n-2}$ 
under a rotation and translation, such that its face of area one is parallel to 
$\pi\in\mathrm{Gr}(n,2)$. Define $S_\pi$ to be the central unit 
square of the plate $P_\pi$.
   When $P_\pi\cap P_{\pi'}\neq \emptyset$ and 
$\theta(\pi,\pi')=\phi$, we say the plates {\it intersect at major angle~$\phi$.} 
Similarly, when $P_\pi\cap P_{\pi'}\neq \emptyset$ and $\un(\pi,\pi')=\phi$,  we 
say the plates {\it intersect at minor angle~$\phi$.}

We will require the following lemmas. The first is due to Mitsis \cite{mi} and is a natural 
extension of an observation of A.~C\'ordoba~\cite{co}. We include the proof for convenience. The third is a natural extension 
of an observation of Wolff~\cite{wo}. 
\begin{lemma}\label{miti} Let $\pi_1, 
{\pi_2}\in\mathrm{Gr}(n,2)$. Then $$|P_{\pi_1}\cap P_{\pi_2}|\ls 
\frac{\de^{n-1}}{\theta(\pi_1,\pi_2)}.$$
\end{lemma} 
\begin{proof} We have that $P_{\pi_1}\cap P_{\pi_2}$ is contained in  
$(R_1\cap R_2)\times R_3$, where $R_1$ and $R_2$ are 2--dimensional rectangles of dimension 
$1\times\de$,  and $R_3$ is an 
$(n-2)$--dimensional cuboid of dimension $1\times\underbrace{\de\times\cdots\times\de}_{n-3}.$
 The rectangles $R_1$ and $ R_2$ intersect at an angle~$\gs \theta(\pi_1,\pi_2)$, so by elementary geometry,
$$
\mathcal{L}^2(R_1\cap R_2)\ls \frac{\de^2}{\theta(\pi_1,\pi_2)},
$$
and the lemma follows.
\end{proof}

\begin{lemma}\label{miti2} Let $\pi_1, {\pi_2}\in\mathrm{Gr}(n,2)$. 
Then $$
|P_{\pi_1}\cap P_{\pi_2}|\ls 
\frac{\de^{n}}{\theta(\pi_1,\pi_2)\underline{\theta}(\pi_1,\pi_2)}.$$
\end{lemma}  
\begin{proof}

By translation we can suppose that the origin is contained in $P_{\pi_1}\cap P_{\pi_2}$. Let 
 $x\in P_{\pi_1}\cap P_{\pi_2}$ and define $l\in\mathrm{Gr}(n,1)$ to be the 
line that passes through the origin and $x$. Define $l_1, l_2\in\mathrm{Gr}(n,1)$ to be the orthogonal projections
of $l$ onto $\pi_1$ and $\pi_2$ respectively.
Now by elementary geometry,
$$
\|\mathrm{proj}_{l}-\mathrm{proj}_{l_1}\|\ls \de/|x|,
$$
and
$$
\|\mathrm{proj}_{l}-\mathrm{proj}_{l_2}\|\ls \de/|x|.
$$
By the triangle inequality,
$$
\underline{\theta}(\pi_1,\pi_2)\le \|\mathrm{proj}_{l_1}-\mathrm{proj}_{l_2}\|\ls \de/|x|,
$$ 
so that
$$
|x|\ls \de/\underline{\theta}(\pi_1,\pi_2).
$$
If we denote the diameter of $P_{\pi_1}\cap P_{\pi_2}$ by $\alpha,$ then we see that $\alpha\ls \de/\underline{\theta}(\pi_1,\pi_2)$.

 As in the previous proof, $P_{\pi_1}\cap P_{\pi_2}$ is contained in  
$(R_1\cap R_2)\times R_3$, where $R_1$ and $R_2$ are 2--dimensional rectangles of dimension 
$1\times\de$,  and $R_3$ is an 
$(n-2)$--dimensional cuboid of dimension $\alpha\times\underbrace{\de\times\cdots\times\de}_{n-3}.$ 
We are able to reduce the length of the long side of $R_3$ as we have a bound 
on the diameter of $P_{\pi_1}\cap P_{\pi_2}$.

The rectangles $R_1$ and $R_2$ intersect at an angle~$\gs \theta(\pi,\pi')$, so by elementary geometry,
$$
\mathcal{L}^2(R_1\cap R_2)\ls \frac{\de^2}{\theta(\pi,\pi')},
$$
and the lemma follows.
\end{proof}

Finally we prove a quantitative version of the fact that three planes 
 intersecting in distinct lines are contained in a 3--plane.

\begin{lemma}\label{trouble} 
 Let $\pi_0,\pi_1,\pi_2\in \mathrm{Gr}(n,2),$ and define 
 $\Sigma=\{x\in\R^n : d(x, S_{\pi_0})> \nu\},$ where $\de<\nu<1$. Suppose that  $P_{\pi_1}, 
P_{\pi_2}$   intersect $P_{\pi_0}$ at major angles $\sim 1$ and minor angles $< \phi$, where 
$\de\le \phi \le 1,$ and 
suppose that
$$P_{\pi_1}\cap P_{\pi_2}\cap \Sigma\neq\emptyset.$$  
 Then there is a 3--space $\la,$ chosen independently of $\pi_2$,  such that
$$
\min_{\pi\subset\la} \theta(\pi,\pi_1)<\phi
\,\,\,\,\,\,\,\,
\mathrm{and}\,\,\,\,\,\,\,\,
\min_{\pi\subset\la} \theta(\pi,\pi_2)\ls\phi/\nu.
$$
 \end{lemma}
 
\begin{proof} By translation we can suppose that the origin is contained in the set
$$
\{ x : d(x,P_{\pi_0}\cap P_{\pi_1})<\de\}\cap S_{\pi_0},
$$
where $d(x,A)=\inf_{y\in A} |x-y|$.
 Now as $\underline{\theta}(\pi_0,\pi_1)<\phi$, there are lines $l_0,l_1\in \mathrm{Gr}(n,1)$ that are contained in      
 $\pi_0$  and $\pi_1$ respectively, such that 
$$
\|\mathrm{proj}_{l_0}-\mathrm{proj}_{l_1}\|<\phi.
$$
Let $l_1'\in\mathrm{Gr}(n,1)$  be the line contained in $\pi_1$ that is orthogonal to $l_1$, and define $\la$ to be the 3--space 
spanned by $\pi_0$ and $l_1'$.
  
  Define 
 $\pi_1'\in\mathrm{Gr}(n,2)$ to be the  plane spanned by $l_0$ and $l_1'$, so that 
 $\pi_1'$ is contained in $\la$. Now by elementary geometry, we have $\theta(\pi_1,\pi_1')<\phi$ and
 $$
\min_{\pi\subset\la} \theta(\pi,\pi_1)<\phi.
$$
  
 Define $T:\R^n\to\R^n$ to be a translation that maps a point in
 $$\{ x : d(x,P_{\pi_0}\cap P_{\pi_2})<\de\}\cap S_{\pi_0}$$ to the origin. We note that $\la$ and $\Sigma$ 
 are essentially invariant under the action of $T,$ so that if $\zeta\in P_{\pi_1}\cap P_{\pi_2}\cap \Sigma$, then there exists a $\zeta'\in \la\cap\Sigma$ such that
$|T(\zeta)-\zeta'|\ls\phi$.

 Now as $\underline{\theta}(\pi_0,\pi_2)<\phi$, there are lines $l_0',l_2\in \mathrm{Gr}(n,1)$
 that are contained in $\pi_0$ and  $\pi_2$ respectively, such that 
$$
\|\mathrm{proj}_{l_0'}-\mathrm{proj}_{l_2}\|<\phi.
$$
Define $\pi_2'\in\mathrm{Gr}(n,2)$ 
 to be the plane that contains $l_0'$ and $\zeta',$ so that $\pi_2'$ is 
 contained in $\la$. Now by elementary geometry, we have
  $\theta(\pi_2,\pi_2')\ls\phi/\nu$ and   
$$
\min_{\pi\subset\la} \theta(\pi,\pi_2)\ls\phi/\nu,
$$
as required.
\end{proof}

\section{The main argument}\label{oneandahalf}

 The {\it (concave) triangle inequality} states that when $p\le1$,  
$$
\|\sum_{k=0}^N f_k\|_p^p \le\sum_{k=0}^{N} \|f_k\|_p^p,
$$
and this will frequently take the following form:
If $\|f_k\|_p^p\ls C$ for all $k$, then 
$$
\|\sum_{k=0}^{\lfloor\log_2(1/\de)\rfloor} f_k\|_p^p\ls C.
$$
Similarly, the {\it pigeonhole principle} will often take the form: If $$
\sum_{k=0}^{\lfloor\log_2(1/\de)\rfloor} \|f_k\|\gs C,
$$ 
then for some $k,$ we have $\|f_{k}\| \gs C.$

The following lemma will be key to the proof of Theorem~\ref{dimension0}.

\begin{lemma}\label{main0}
Let $\Pi$ be a $\de$--separated subset of a  1--dimensional
 manifold 
$\M$ in $\mathrm{Gr}(3,2)$. Then
\begin{equation*}
\|\sum_{\pi\in\Pi}\chi_{P_\pi}\|_{2}\ls 1.
\end{equation*}
\end{lemma}

\begin{proof} 
 We note that as
$$
 \|\sum_{\pi\in\Pi}\chi_{P_\pi}\|_{2}^2=\|\sum_{\pi\in\Pi}\sum_{\pi'\in\Pi}\chi_
{P_\pi}\chi_{P_{\pi'}}\|_{1}, $$
  it will suffice to show
 $$ 
\|\sum_{\pi\in\Pi}\sum_{\pi'\in\Pi}\chi_{P_\pi}\chi_{P_{\pi'}}\|_{1}\ls 1.  
$$
Now 
 $$
  \sum_{\pi\in\Pi}\sum_{\pi'\in\Pi}\chi_{P_\pi}\chi_{P_{\pi'}}\le
  \sum_{k=0}^{\lfloor\log_2(1/\de)\rfloor}\sum_{\pi,\pi':\theta(\pi,\pi')\sim2^{-k
}}\chi_{P_\pi}\chi_{P_{\pi'}}+\sum_{\pi\in\Pi}\chi_{P_\pi},  $$  as $\Pi$ is 
$\de$--separated. Thus, by   the triangle inequality, it will suffice to show 
  \begin{equation}\label{dis0}
  \|\sum_{\pi,\pi':\theta(\pi,\pi')\sim2^{-k}}\chi_{P_\pi}\chi_{P_{\pi'}}\|_{1}\ls 1  
  \end{equation}
   for all $k$ less than or equal to 
$\lfloor\log_2(1/\de)\rfloor$, and   
\begin{equation*}
  \|\sum_{\pi\in\Pi}\chi_{P_\pi}\|_{1}\ls 1.
  \end{equation*}
  Again, by the triangle inequality, 
  $$
   \|\sum_{\pi\in\Pi}\chi_{P_\pi}\|_{1}\le \sum_{\pi\in\Pi}|P_\pi|\ls 1,
   $$
    as $\#\Pi\ls \de^{-1}$ and $|P_\pi|=\de.$ Thus it 
remains  to show (\ref{dis0}) for each $k$, which we now consider to be fixed.  
  
  Using the metric $\theta$ on the 1--dimensional manifold $\mathcal{M}$,  we can 
  cover $\Pi$ by a constant multiple of $2^k$  balls $\{B_j\}$, with radius a constant 
  multiple of $2^{-k}.$  We can also choose the cover 
so that if $\theta(\pi,\pi')\sim2^{-k}$, then  $\pi$ and $\pi'$ are both 
contained in some $B_j$. Hence, by the triangle inequality, it will suffice to 
prove    
\begin{equation}\label{roar0}  \|\sum_{\substack{\pi,\pi'\in 
B_j\cap\Pi\\  
\theta(\pi,\pi')\sim2^{-k}}}\chi_{P_\pi}\chi_{P_{\pi'}}\|_{1}\ls 
2^{-k} \end{equation}
  for all $k$ less than or equal to 
$\lfloor\log_2(1/\de)\rfloor$, and each ball $B_j$.     
  
  Without loss of generality, we can suppose that $B_j$ is centered on the $x_1x_2$--plane. Define the 
dilation $L:\R^3\to \R^3$ by  $$
  L(x_1,x_2,x_3)=(x_1,x_2,2^kx_3).  
  $$
  We scale our geometric configuration by $L$, so that 
  \begin{align}\label{meek0}
  \|\sum_{\substack{\pi,\pi'\in B_j\cap\Pi\\
  \theta(\pi,\pi')\sim2^{-k}}}\chi_{P_\pi}\chi_{P_{\pi'}}\|_{1}&=2^{-k}
\|\sum_{\substack{\pi,\pi'\in B_j\cap\Pi\\  
\theta(\pi,\pi')\sim2^{-k}}}\chi_{L(P_\pi)}\chi_{L(P_{\pi'})}\|_{1}.
  \end{align}  
 Now as we have essentially 
  changed $\de$  to $2^k\de$, and $\theta(L(\pi),L(\pi'))\sim 1,$  
  if we can prove (\ref{roar0}) 
  when $k=0$, then 
  $$
  \|\sum_{\substack{\pi,\pi'\in B_j\cap\Pi\\  
\theta(\pi,\pi')\sim2^{-k}}}\chi_{L(P_\pi)}\chi_{L(P_{\pi'})}\|_{1}
\le 
C_\epsilon(2^k\de)^{-\epsilon}< C_\epsilon\de^{-\epsilon}\ls 1, 
  $$
   so that by (\ref{meek0}), 
  $$
  \|\sum_{\substack{\pi,\pi'\in B_j\cap\Pi\\
  \theta(\pi,\pi')\sim2^{-k}}}\chi_{P_\pi}\chi_{P_{\pi'}}\|_{1}\ls 2^{-k}.
  $$
  Thus   
  it will suffice to prove (\ref{roar0}) when $k=0$.
  
  Now 
$$
 \|\sum_{ \theta(\pi,\pi')\sim 1}\chi_{P_\pi}\chi_{P_{\pi'}}\|_{1}\le 
\sum_{ \theta(\pi,\pi')\sim 1}|P_\pi\cap P_{\pi'}|\ls 
(\#\Pi)^2\de^{2}\ls1, $$ 
by Lemma~\ref{miti}, and we are done. 
\end{proof}

The following lemma will be key to the proof of Theorem~\ref{dimension}.

\begin{lemma}\label{main}
Let $\Pi$ be a $\de$--separated subset of a  2--dimensional
curved manifold
$\M$ in $\mathrm{Gr}(4,2)$. Then
\begin{equation*}\label{mainy}
\|\sum_{\pi\in\Pi}\chi_{P_\pi}\|_{5/3}\ls \de^{-1/5}.
\end{equation*}
\end{lemma}

\begin{proof} 
The proof is based on the ideas of Wolff \cite{wo}.
The key geometric fact used there, is that three lines intersecting in distinct points are contained in a plane.
The corresponding fact here, is that three planes intersecting in distinct lines are contained in a 
3--plane. Unfortunately, the intersection between two planes can be a point as well as a line. 

In order to deal with the different types of intersection, 
we use the bilinear reduction of T. Tao, A. Vargas and  L. Vega~\cite{ta2}, 
which can also be found in~\cite{ta1}.  This enables us to quantify, using the minor angle, how
near the planes are to intersecting in lines. When the planes are not intersecting in lines,
we are able to use  Lemma~\ref{miti2}, in place of Lemma~\ref{miti}, in compensation. 

We make the bilinear reduction. Essentially this means we will begin by attempting to
copy the proof of Lemma~\ref{main0}. We note that as
$$
 \|\sum_{\pi\in\Pi}\chi_{P_\pi}\|_{5/3}^{5/3}=\|\sum_{\pi\in\Pi}\sum_{\pi'\in\Pi}\chi_
{P_\pi}\chi_{P_{\pi'}}\|_{5/6}^{5/6}, $$
  it will suffice to show
 $$ 
\|\sum_{\pi\in\Pi}\sum_{\pi'\in\Pi}\chi_{P_\pi}\chi_{P_{\pi'}}\|_{5/6}^{5/6}\ls 
\de^{-1/3}.  
$$   Now 
 $$
  \sum_{\pi\in\Pi}\sum_{\pi'\in\Pi}\chi_{P_\pi}\chi_{P_{\pi'}}\le
  \sum_{k=0}^{\lfloor\log_2(1/\de)\rfloor}\sum_{\pi,\pi':\theta(\pi,\pi')\sim2^{-k
}}\chi_{P_\pi}\chi_{P_{\pi'}}+\sum_{\pi\in\Pi}\chi_{P_\pi},  $$  as $\Pi$ is 
$\de$--separated. Thus, by   the triangle inequality, it will suffice to show  
  \begin{equation}\label{dis}
  \|\sum_{\pi,\pi':\theta(\pi,\pi')\sim2^{-k}}\chi_{P_\pi}\chi_{P_{\pi'}}\|_{5/6}^
{5/6}\ls \de^{-1/3}  \end{equation}
   for all $k$ less than or equal to 
$\lfloor\log_2(1/\de)\rfloor$, and   
\begin{equation*}\label{dis2}
  \|\sum_{\pi\in\Pi}\chi_{P_\pi}\|_{5/6}^{5/6}\ls \de^{-1/3}.
  \end{equation*}
   Again by the triangle inequality,
  $$
   \|\sum_{\pi\in\Pi}\chi_{P_\pi}\|_{5/6}^{5/6}\le \|\sum_{\pi\in\Pi}\chi_{P_\pi}\|_{1}\le \sum_{\pi\in\Pi}|P_\pi|\ls 1<\de^{-1/3},
   $$
  as $\#\Pi\ls\de^{-2}$ and $|P_\pi|=\de^{2}$, so it 
remains  to show (\ref{dis}) for each $k$, which we now consider to be fixed.  
 
Using the metric $\theta$ on the 2--dimensional manifold $\mathcal{M}$,  we can 
  cover $\Pi$ by a constant multiple of $2^{2k}$  balls $\{B_j\}$, with radius a constant 
  multiple of $2^{-k}.$  We can also choose the cover 
so that if $\theta(\pi,\pi')\le2^{-k}$, then  $\pi$ and $\pi'$ are both 
contained in some $B_j$.
  Hence, by the triangle inequality, it will suffice to 
prove    
\begin{equation}\label{roar}  \|\sum_{\substack{\pi,\pi'\in 
B_j\cap\Pi\\  
\theta(\pi,\pi')\sim2^{-k}}}\chi_{P_\pi}\chi_{P_{\pi'}}\|_{5/6}^{5/6}\ls 
2^{-2k}\de^{-1/3}  \end{equation}
  for all $k$ less than or equal to 
$\lfloor\log_2(1/\de)\rfloor$, and each ball $B_j$.   
     
  Without loss of generality, we can suppose that $B_j$ is centered on the $x_1x_2$--plane. Define the 
dilation $L:\R^4\to \R^4$ by  $$
  L(x_1,x_2,x_3,x_4)=(x_1,x_2,2^kx_3,2^kx_4).  
  $$
  As $\det L=2^{2k}$, if we scale our geometric configuration by $L$, then 
  \begin{align}\label{belle}
  \|\sum_{\substack{\pi,\pi'\in B_j\cap\Pi\\
  \theta(\pi,\pi')\sim2^{-k}}}\chi_{P_\pi}\chi_{P_{\pi'}}\|_{5/6}^{5/6}&=2^{-2k}
\|\sum_{\substack{\pi,\pi'\in B_j\cap\Pi\\  
\theta(\pi,\pi')\sim2^{-k}}}\chi_{L(P_\pi)}\chi_{L(P_{\pi'})}\|_{5/6}^{5/6}.
  \end{align}  
  
 Essentially we have 
  changed $\de$  to $2^k\de,$ and $\theta(L(\pi),L(\pi'))\sim 1,$ so that 
  if we can prove (\ref{roar}) 
  when $k=0$, then 
  \begin{equation}\label{sebastian}
  \|\sum_{\substack{\pi,\pi'\in B_j\cap\Pi\\  
\theta(\pi,\pi')\sim2^{-k}}}\chi_{L(P_\pi)}\chi_{L(P_{\pi'})}\|_{5/6}^{5/6}
\le 
C_\epsilon(2^k\de)^{-1/3-\epsilon}< C_\epsilon\de^{-1/3-\epsilon}\ls \de^{-1/3}. 
  \end{equation}
   Now by combining (\ref{belle}) and (\ref{sebastian}), 
  $$
  \|\sum_{\substack{\pi,\pi'\in B_j\cap\Pi\\
  \theta(\pi,\pi')\sim2^{-k}}}\chi_{P_\pi}\chi_{P_{\pi'}}\|_{5/6}\ls 2^{-2k}\de^{-1/3},
  $$
  so 
  it will suffice to prove (\ref{roar}) when $k=0$.

  Now as   
$$  \sum_{\pi,\pi'}\chi_{P_\pi}\chi_{P_{\pi'}}\le  
\sum_{k=0}^{\lfloor\log_2(1/\de)\rfloor}\sum_{\pi,\pi':\un(\pi,\pi')\sim2^{-k}}
\chi_{P_\pi}\chi_{P_{\pi'}}+\sum_{\pi,\pi':\un(\pi,\pi')< 
\de}\chi_{P_\pi}\chi_{P_{\pi'}},  $$ 
 by the triangle inequality, it will suffice 
to show  
\begin{equation}\label{points}  \|\sum_{\substack{\theta(\pi,\pi')\sim1\\ \un(\pi,\pi')\sim 
2^{-k}}}\chi_{P_\pi}\chi_{P_{\pi'}}\|_{5/6}^{5/6}\ls \de^{-1/3}  
\end{equation}
for all $k$ less than or equal to $\lfloor\log_2(1/\de)\rfloor,$ 
and   
\begin{equation}\label{lines}
  \|\sum_{\substack{\theta(\pi,\pi')\sim 1\\ \un(\pi,\pi')< \de}}\chi_{P_\pi}\chi_{P_{\pi'}}\|_{5/6}^{5/6}\ls 
\de^{-1/3}.  
\end{equation}
  In (\ref{points}) the planes are only intersecting  in points, 
and  in (\ref{lines}) the planes are almost intersecting in lines.  
  
   To prove (\ref{points}), we fix $k$ and 
   define $F$ by  
    $$
  F=\{x\in \R^4:\sum_{\pi,\pi'\in\Pi:(*)}\chi_{P_\pi}\chi_{P_{\pi'}}(x)\sim \lambda\},
  $$
  where $(*)$ denotes the conditions $\theta(\pi,\pi')\sim 1$ and $\un(\pi,\pi')\sim 2^{-k}$. We have 
   $$ 
\|\sum_{\pi,\pi'\in\Pi:(*)}\chi_{P_\pi}\chi_{P_{\pi'}}\|_{5/6}^{5/6}
\ls \sum_{\lambda} \lambda^{5/6}|F|,  $$ where 
 $\lambda$ ranges dyadically up to a constant multiple of $\de^{-4}$. Thus, by 
 the triangle inequality, 
  it will suffice to show the weak type inequality
\begin{equation}\label{res}
  \lambda^{5/6}|F|\ls \de^{-1/3}.
  \end{equation}
  We can  assume that $|F|$ is greater than $\de^3,$ as otherwise we are 
done.
  
  Define $\Pi_\lb$  by
  $$
  \Pi_\lb=\{\pi\in\Pi: |P_\pi\cap F |\sim \lb| P_\pi|\}.
  $$
  We will use the pigeonhole principle to find a single plate that intersects 
many other plates that have density $\lb\gs \sqrt{\lambda} |F|.$  
  
    By definition,   $$
  \int_{F}\sum_{\pi, \pi'\in\Pi:(*)}\chi_{P_\pi}\chi_{P_{\pi'}} \gs \lambda|F|,  $$
 so that 
  $$
  \sum_{\lb, \lb'}\sum_{\pi\in\Pi_\lb}\sum_{ \pi'\in\Pi_{\lb'}:(*)} |P_\pi\cap P_{\pi'}\cap 
F|=\int_{F}\sum_{\pi,\pi\in\Pi:(*)}\chi_{P_\pi}\chi_{P_{\pi'}}\gs \lambda|F|\gs \de^3,  $$
  where the sums over $\lb$ and $\lb'$ range dyadically from zero to one. The summands 
where $\lb$ or $\lb'$ is less than a large power of $\de$ can be absorbed by the  larger   
summands. Thus, by the pigeonhole principle, we can find $\lb$ and $\lb'$, which we now 
fix,  such that  
\begin{equation}\label{pig}  
\int_{F}\sum_{\pi\in\Pi_\lb}\sum_{\pi'\in\Pi_{\lb'}:(*)} 
\chi_{P_\pi}\chi_{P_{\pi'}} \gs \lambda|F|.  
\end{equation}   
Let
$$F'=\{x\in F:\sum_{\pi\in\Pi_\lb}\sum_{\pi'\in\Pi_{\lb'}:(*)} 
\chi_{P_\pi}\chi_{P_{\pi'}}\gs\lambda\},$$
 so that by (\ref{pig}) and the pigeonhole principle,
$$
\int_{F'}\sum_{\pi\in\Pi_\lb}\sum_{\pi'\in\Pi_{\lb'}:(*)} 
\chi_{P_\pi}\chi_{P_{\pi'}} \gs \lambda|F|.  
$$
By definition,
  $$
  \sum_{\pi\in\Pi_{\lb}}\sum_{\pi'\in\Pi_{\lb'}:(*)}\chi_{P_\pi}\chi_{P_{\pi'}}\le\lambda
  $$
  on $F$, so that  $|F'|\gs|F|.$ Now as $$
   \sum_{\pi\in\Pi_\lb}\sum_{\pi'\in\Pi_{\lb'}:(*)} 
\chi_{P_\pi}\chi_{P_{\pi'}}\gs\lambda$$
 on $F'$,  
  we can suppose that
  $$
  \sum_{\pi\in\Pi_\lb}\chi_{P_\pi}\gs\sqrt{\lambda}
  $$
  on some $F''\subset F'$, where $|F''|\gs|F'|.$ Thus
  $$
  \int_{F}\sum_{\pi\in\Pi_\lb}\chi_{P_\pi}\ge \int_{F''}\sum_{\pi\in\Pi_\lb}\chi_{P_\pi}
  \gs \sqrt{\lambda}|F''|\gs \sqrt{\lambda}|F|.
  $$
By definition,  $\lb\de^2\gs |P_\pi\cap F|$ for all $\pi\in\Pi_{\lb},$ so that
$$
\sum_{\pi\in\Pi_\lb}\lb\de^2\gs \sum_{\pi\in\Pi_\lb}|P_\pi\cap F|=\int_{F}\sum_{\pi\in\Pi_\lb}\chi_{P_\pi}\gs \sqrt{\lambda}|F|.
$$
Now as $\#\Pi_\lb\ls\de^{-2}$, we have $\lb\gs\sqrt{\lambda}|F|.$
  
  From (\ref{pig}), we also have 
  $$
  \sum_{\pi'\in\Pi_{\lb'}}\int_{P_{\pi'}}\sum_{\pi\in\Pi_{\lb}:(*)} \chi_{P_\pi} \gs 
\lambda|F|,  $$
   so that, there exists a
$\pi_0\in\Pi_{\lb'}$,  such that  
\begin{equation}\label{pigs}
  \int_{P_{\pi_0}}\sum_{\pi\in\Pi_{\lb}:(*)} \chi_{P_\pi} \gs 
\de^{2}\lambda|F|,  
\end{equation}
where $(*)$ denotes the conditions $\theta(\pi,\pi_0)\sim 1$ and $\un(\pi,\pi_0)\sim 2^{-k}$. 
  By Lemma~\ref{miti2},
  \begin{equation}\label{pigss}
  \int_{P_{\pi_0}} \chi_{P_\pi}=|P_\pi\cap P_{\pi_0}|\ls 2^{k}\de^{4},
  \end{equation}
  so that
  if we define $\mathcal{P}$ by
  $$\mathcal{P}=\{\pi\in\Pi_\lb: \theta(\pi,\pi_0)\sim1,\,\,\un(\pi,\pi_0)\sim 2^{-k},\,\,P_\pi\cap P_{\pi_0}\neq 
  \emptyset \},$$
  then by (\ref{pigs}), we have
  \begin{equation*}\label{argh2}
  \#\mathcal{P}\gs 2^{-k}\de^{-2}\lambda|F|.
  \end{equation*}
 Thus we have a large set of planes, with density $\lb\gs \sqrt{\lambda} |F|,$ that intersect $\pi_0$.

As $\theta(\pi,\pi_0)\sim 1$ for all $\pi\in\PP$, we can define $\Sigma=\{x\in\R^4 : d(x,S_{\pi_0})>\lb/C\}$ for 
some sufficiently large constant $C$, so that 
 $$
 \int_{F}\chi_{P_\pi\cap\Sigma}\gs \lb\de^{2}.
 $$
  Summing over $\PP$, we have
  $$
  \int_{F}\sum_{\pi\in\PP}\chi_{P_\pi\cap\Sigma}\gs 
\lb\de^{2}\#\PP,  $$
  so that by the Cauchy--Schwarz inequality,
  \begin{equation}\label{argh}
  \|\sum_{\pi\in\PP}\chi_{P_\pi\cap\Sigma}\|_2\gs\frac{\lb\de^{2}\#\PP}{|F|^{1/2}}.  
\end{equation}  

We will use the geometry of the construction to bound the left hand side of (\ref{argh}) from above,  in 
order to obtain  (\ref{res}). 
  If we fix a $\pi_1\in\PP$, then by Lemma~\ref{trouble}, there 
exists a 3--space $\la$, such that 
\begin{equation}\label{mere}
\min_{\pi\subset\la} \theta(\pi,\pi_1)<2^{-k},
\end{equation}
and if $\pi_2\in \mathcal{P}$ and $P_{\pi_1}\cap P_{\pi_2}\cap \Sigma\neq\emptyset,$ then
\begin{equation}\label{mere2}
\pi_2\in \mathcal{M}^{\lb^{-1}2^{-k}}_\Lambda, 
\end{equation}
where $\mathcal{M}^{\lb^{-1}2^{-k}}_\Lambda=\{\pi'\in\mathcal{M}:\min_{\pi\subset\la}\theta(\pi,\pi')<\lb^{-1}2^{-k}\}.$ 
%\begin{equation}\label{mere2}
%\min_{\pi\subset\la} \theta(\pi,\pi_2)\ls\lb^{-1}2^{-k}.
%\end{equation}

Now, as $\mathcal{M}$ is curved, there exists a constant $C$ such that
\begin{equation*}\label{pope2}
\frac{\mu(\mathcal{M}^{\epsilon_1}_{\la}\cap 
 B_{\epsilon_2}(\pi))}{\epsilon_1\epsilon_2}<C
\end{equation*}
for all $\epsilon_1,\epsilon_2>0$, 3--spaces $\Lambda$, and $\pi\in\mathcal{M}.$
 Thus, as $\PP$ is 
 $\de$--separated, we see that
$$
  \#\{ \pi_2\in\PP :\theta(\pi_1,\pi_2)\le 2^{-j},\,\,\, P_{\pi_1}\cap P_{\pi_2}\cap 
\Sigma\neq\emptyset\}\ls   
\lb^{-1}2^{-k}2^{-j}\de^{-2}, 
$$
where we take $j$ to be less than or equal to $\lfloor\log_2(1/\de)\rfloor$.
By Lemma~\ref{miti}, 
 $$|P_{\pi_1}\cap P_{\pi_2}\cap\Sigma|\ls 2^j\de^{3}$$  
when $\theta(\pi_1,\pi_2)\sim 2^{-j}$, so that
 $$
 \sum_{\pi_2\in\PP: \theta(\pi_1,\pi_2)\sim 2^{-j}}|P_{\pi_1}\cap P_{\pi_2}\cap \Sigma|\ls \lb^{-1}2^{-k}\de.
 $$
 Hence, by the triangle inequality,
 $$
 \sum_{\pi_2\in\PP}|P_{\pi_1}\cap P_{\pi_2}\cap \Sigma|\ls \lb^{-1}2^{-k}\de,
 $$
 so that by summing over $\PP$, we have
 $$
  \sum_{\pi_1\in\PP}\sum_{\pi_2\in\PP}|P_{\pi_1}\cap P_{\pi_2}\cap \Sigma|\ls \lb^{-1}2^{-k}\de\#\PP.
 $$
 Now
 $$
 \| \sum_{\pi\in\PP}\chi_{P_{\pi}\cap \Sigma}\|^2_2=\sum_{\pi_1\in\PP}\sum_{\pi_2\in\PP}|P_{\pi_1}\cap P_{\pi_2}\cap \Sigma|,
 $$
 so that
 \begin{equation}\label{second}
 \| \sum_{\pi\in\PP}\chi_{P_{\pi}\cap \Sigma}\|_2\ls (\lb^{-1}2^{-k}\de\#\PP)^{1/2}.
 \end{equation}
 
  Combining the equations (\ref{argh}) and (\ref{second}),  
   and the fact that $\#\mathcal{P}\gs2^{-k}\de^{-2}\lambda|F|$, we obtain
 $$
 \lb^3\lambda \ls \de^{-1}.
 $$
  Using the fact that $\nu\gs \sqrt{\lambda}|F|$, we have
  $$
  \lambda^{5/2}|F|^{3}\ls \de^{-1},
  $$
  and we take the third root to obtain
   (\ref{res}) as required.
   
   To prove (\ref{lines}), we argue in the same way. We let 
    $(*)$ denote the conditions $\theta(\pi,\pi')\sim 1$ and $\un(\pi,\pi')<\de$.
   We apply Lemma~\ref{miti} in place of Lemma~\ref{miti2}, so that the estimate (\ref{pigss}) becomes  
  \begin{equation*}
  \int_{P_{\pi_0}} \chi_{P_\pi}=|P_\pi\cap P_{\pi_0}|\ls \de^{3},
  \end{equation*}
  and $\#\mathcal{P}\gs \de^{-1}\lambda|F|.$ The expressiones in (\ref{mere}) and (\ref{mere2}) are changed to
  \begin{equation*}
\min_{\pi\subset\la} \theta(\pi,\pi_1)<\de
\,\,\,\,\,\mathrm{ and }\,\,\,\,\,
\pi_2\in \mathcal{M}^{\lb^{-1}\de}_\Lambda,
\end{equation*}
  so that
  $$
  \| \sum_{\pi\in\PP}\chi_{P_{\pi}\cap \Sigma}\|_2\ls (\lb^{-1}\de^{2}\#\PP)^{1/2}.
  $$
  As before, we combine these equations with (\ref{argh}), which is unchanged,  so that
   $$
  \lambda^{5/2}|F|^{3}\ls \de^{-1},
  $$
    and we are done.
\end{proof}

\section{Proof of Theorems \ref{dimension0} and \ref{dimension}}\label{two}

The following argument is well known, and can be found in \cite{bo}. Let $\{B(x_j,r_j)\}$ 
be a covering of a $\mathcal{K}^n_2$-set $E$, where $r_j\le1/4$. We are required to show 
$$
\sum_j r_j^{\frac{n+3}{2}-\epsilon}\ge C_\epsilon>0
$$
for all $\epsilon>0$.
 
Define $E_k$ by
$$
E_k=E\cap\bigcup_{r_j\sim2^{-k}}B(x_j,r_j)
$$
for $k\ge 2$.
As $E$ is a $\mathcal{K}^n_2$-set, for all $\pi\in\mathcal{M}$ there is a square 
$S_\pi\subset E$ with a corresponding plate $P_\pi$.
By the pigeonhole principle, for all $\pi\in\mathcal{M}$   
there is a $k_\pi$ such that
$$
\mathcal{L}^2(S_\pi\cap E_{k_\pi})\ge \frac{1}{k_{\pi}^2}\gs 1.
$$
Thus $\mathcal{M}=\bigcup_{k\ge 2}\mathcal{M}_k,$ where 
$$
\mathcal{M}_k=\{\pi\in\mathcal{M}:\mathcal{L}^2(S_\pi\cap E_k)\ge \frac{1}{k^2}\}.
$$
By the pigeonhole principle again, there exists a $k_0$ such that
$$
\mu(\mathcal{M}_{k_0})\ge \frac{1}{k_0^2}\gs 1,
$$
where $\mu$ is the induced measure on $\mathcal{M}$, normalised to have total mass 1.  

Let $\delta=2^{-k_0},$ and let $\Pi$ be a maximal $\de$--separated 
subset of $\mathcal{M}_{k_0},$ so that $\#\Pi\gs \de^{2-n}$.
 Define
$
J=\{ j : r_j\sim \de\}
$
and $E_\de=\bigcup_{j\in J}B(x_j,2\de),$ so that 
$$
\int_{E_\de}\chi_{P_\pi}\gs\de^{n-2}\mathcal{L}^2(S_\pi\cap E_{k_0})\ge\frac{\de^{n-2}}{k_0^2}\gs\de^{n-2}.
$$
Now as  $\#\Pi\gs \de^{2-n}$, we have
$$
\sum_{\pi\in \Pi} \int_{E_\de}\chi_{P_\pi}\gs 1,
$$
 so that by H\"{o}lder's inequality,
$$
|E_\de|^{\frac{2}{n+1}}\|\sum_{\pi\in \Pi}\chi_{P_\pi}\|_{\frac{n+1}{n-1}}\gs 1.
$$
By Lemmas~\ref{main0}~and~\ref{main}, we have
$$
\|\sum_{\pi\in \Pi}\chi_{P_\pi}\|_{\frac{n+1}{n-1}}\ls \de^{\frac{3-n}{n+1}},
$$
where $n=3$ or $4$, so that
$$
|E_\de|\gs \de^{\frac{n-3}{2}}.
$$

On the other hand, we have $\#J\de^n\gs|E_\de|$, so that $$\#J\ge C_\epsilon\de^{-\frac{n+3}{2}+\epsilon}$$
for all $\epsilon>0$. Hence, when $n=3$ or $4$, 
 $$
 \sum_j r_j^{\frac{n+3}{2}-\epsilon}\ge \#J \left(\frac{\de}{2}\right)^{\frac{n+3}{2}-\epsilon}\ge C_\epsilon'>0
 $$ 
for all $\epsilon>0$, and we are done. \hfill $\square$

\section{Sharpness in the complex case}\label{three}

We construct an example, inspired by the Heisenberg group example in \cite{ka}, 
to show that Theorem~\ref{cork} is sharp. Define $E\subset\C^4$ 
by
$$
E=\{(z_1,z_2,z_3,z_4)\in\C^4 : 
\mathrm{Im}(z_1\overline{z_2})=\mathrm{Im}(z_3\overline{z_4})\},$$
so that $\dim_H(E)=7$. 
Define the planes $\pi_{u,v}$ by
$$
\pi_{u,v}=
\{(1,-|v|\mathrm{Im}(u)i,u,|v|)z+(1,|v|\mathrm{Im}(u)i,\frac{\overline{u}v}
{|v|},v)z':z,z'\in\C\},
$$
and the manifold $\mathcal{M}$ by
$$\mathcal{M}=\{\pi_{u,v}: \,u,v\in\C, \,\,\mathrm{Im}(u)\neq0,\,\,\mathrm{Im}(v)\neq0\}.
$$
It is not hard to calculate that the planes are contained in $E$, 
and it is clear that $\mathcal{M}$ is a 
2--dimensional manifold as a subset of $\mathrm{Gr}(4,2)$ over~$\C$. 

It remains to 
show that $\mathcal{M}$ is curved. It will suffice to show that $\{\pi\in\mathcal{M} : 
\pi\subset\Lambda\}$ is no more than 1--dimensional for all 3--spaces $\Lambda$. 
Now if a plane in $\mathcal{M}$ is contained in a 3--space $\Lambda,$ then 
\begin{align*} 
(1,-|v|\mathrm{Im}(u)i,u,|v|)\cdot(a_1,a_2,a_3,a_4)=a_1-a_2|v|
\mathrm{Im}(u)i+a_3u+a_4|v|=0 
\end{align*}
and
$$ 
(1,|v|\mathrm{Im}(u)i,\frac{\overline{u}v}{|v|},v)\cdot(a_1,a_2,a_3,a_4)=a_1+a_2
|v|\mathrm{Im}(u)i+a_3\frac{\overline{u}v}{|v|}+a_4v=0$$
for some normal $(a_1,a_2,a_3,a_4)\in\C^4$. We can multiply the second equation 
by $|v|/v$ and subtract it from the first to solve for $\mathrm{Im}(u)$ in terms 
of $v$. Substituting back into the first equation, we fix $u$ in terms of $v$, so 
that the set of planes contained in $\Lambda$ is parametrized by a single 
variable. Thus the restriction of $\mathcal{M}$ to a 3--space is no more than 1--dimensional, 
and $E$ is a $\mathcal{K}^4_2$-set.

Thus, in the complex case, the curvature 
condition is not sufficient to guarantee nontrivial intersection, even before 
translation. This example does not extend to the reals (or the finite fields, as 
there is no square root), and the curvature condition is stronger over $\R$ than over $\C$.
It seems possible that the real and complex cases are qualitatively 
different. Thus the problem of sharp lower bounds for the Hausdorff dimension 
of real $\mathcal{K}^n_2$-sets is open and interesting for~$n\ge 4.$

\section{Final remarks}

We could define the $\mathcal{K}^n_{2}$-sets so that they only 
contain a unit square parallel to each direction plane, and not necessarily the whole plane.

Theorems~\ref{dimension0} and~\ref{dimension} would extend in a natural way to 
$\mathcal{K}^n_{k}$-sets, where $k=n-1$ or $n-2$, and the planes are replaced by $k$--planes. 
We note that the $\mathcal{K}^4_{2}$-set over $\C$ of the previous section may not be a 
$\mathcal{K}^8_{4}$-set over $\R$, as the curvature condition is stronger over $\R$. 
   
%We could perhaps define the
% $\mathcal{K}^n_{k}$-sets with more general direction sets than manifolds. We could take 
% $\mathcal{M}$ to be a curved subset of $\mathrm{Gr}(n,k)$ with Minkowski dimension $n-k.$
%  Curved in this context would mean that the Minkowski dimension of
 %$\{\pi\in\mathcal{M}:\pi\subset\la\}$ is less than or equal to $m-k$ for all $m$--spaces $\la$.
   
   It also 
 seems likely that we could adapt the proofs to obtain the corresponding maximal function estimates. 
 We neglect
 these potential generalizations mainly for expository purposes.
  \vspace{1em}

\thanks{Many thanks to Gerd Mockenhaupt for suggesting a similar planar variant problem 
which led to the consideration of this one. 
Thanks to Peter Sj\"ogren, Thomas Duyckaerts, Fulvio Ricci, Ana Vargas and Toby Bailey for 
helpful conversations, 
and thanks to an anonymous referee for pointing out a mistake in an earlier version. 
%Thanks also to the anonymous referee(s) for helpful comments. 
}

\end{document}